\documentclass[a4paper]{article}
\usepackage{amsmath}

\newtheorem{theorem}{Theorem}

\numberwithin{theorem}{section}
\numberwithin{equation}{section}
\numberwithin{table}{section}

\begin{document}

\title{\textbf{Spinor Darboux Equations of Curves in Euclidean 3-Space}}
\author{\.{I}lim Ki\c{s}i$^1$, Murat Tosun$^2$ \\
}

\date{}
\maketitle

\small
$^1$Department of Mathematics, Kocaeli University, TR41380 Kocaeli, Turkey\\
$^2$Department of Mathematics, Sakarya University, TR54187 Sakarya, Turkey

\begin{abstract}
In this paper, the spinor formulation of Darboux frame on an
oriented surface is given. Also, the relation between the spinor
formulation of Frenet frame and Darboux frame are obtained.
\end{abstract}




\section{Introduction }
Spinors in general were discovered by Elie Cartan in 1913 \cite{c}. Later,
spinors were adopted by quantum mechanics in order to study the
properties of the intrinsic angular momentum of the electron and
other fermions. Today, spinors enjoy a wide range of physics
applications. In mathematics particularly in differential geometry
and global analysis, spinors have since found  broad applications
to algebraic and differential topology, symplectic geometry, gauge
theory, complex algebraic geometry, index theory \cite{b,n}.

In the differential geometry of surfaces a Darboux frame is a
naturel moving frame constructed on a surface. It is the analog of
the Frenet- Serret frame as applied to surface geometry. A Darboux
frame exists at any non-umbilic point of a surface embedded in
Euclidean space. The construction of Darboux frames on the
oriented surface first considers frames moving along a curve in
the surface, and then specializes when the curves move in the direction of
the principal curvatures \cite{m}.

In \cite{c,cb}, the triads of mutually orthogonal unit vectors were
expressed in terms of a single vector with two complex components,
called a spinor. In the light of the existing studies, the Frenet
equations reduce to a single spinor equation, equivalent to the
three usual vector equations, is a consequence of the relationship
between spinors and orthogonal triads of vectors. The aim of this
paper is to show that the Darboux equations can be expressed a
single equation for a vector with two complex components.

\section{Preliminaries}
The Euclidean 3-space   provided with standard flat metric given
by
\begin{equation*}
\langle,\rangle=dx_1^2+dx_2^2+dx_3^2
\end{equation*}
where $(x_1,x_2,x_3)$  is a rectangular coordinate system of $E^3$.
Recall that, the norm of an arbitrary vector $x\in{E^3}$ is given
by $\|x\|=\sqrt{\langle{x,x}\rangle}$. Let $\alpha$  be a curve in
Euclidean 3- space. The curve $\alpha$  is called a unit speed
curve if velocity vector $\alpha'$  of  $\alpha$ satisfies
$\|\alpha'\|=1$.

Let us denote by $T(s)$,$N(s)$ and $B(s)$ the unit tangent vector, unit normal vector and unit binormal vector of $\alpha$ respectively. The Frenet Trihedron is the collection of $T(s)$,$N(s)$ and $B(s)$. Thus, the Frenet formulas are as
\begin{align*}
\frac{dT}{ds}&=\kappa(s)N(s) \\
\frac{dN}{ds}&=-\kappa(s)T(s)+\tau(s)B(s) \\
\frac{dB}{ds}&=-\tau(s)N(s)
\end{align*}
Here, the curvature is defined to be  $\kappa(s)=\|T'(s)\|$ and the torsion is the function $\tau$ such that $B'(s)=-\tau(s)N(s)$ \cite{dc}.

The group of rotation about the origin denoted by $SO(3)$ in $R^3$ is homomorphic to the group of unitary complex $2\times2$ matrices with unit determinant denoted by $SO(2)$. Thus, there exits a two-to-one homomorphism of $SO(2)$ onto $SO(3)$. Whereas the elements of $SO(3)$ act the vectors with three real components, the elements of $SO(2)$ act on vectors with two complex components which are called spinors \cite{sw,g}.
 In this case, we can define a spinor
\begin{equation}  \label{eq:twoone}
\psi=\left(
       \begin{array}{c}
         \psi_1 \\
         \psi_2 \\
       \end{array}
     \right)
\end{equation}
by means of three vectors $a,b,c\in{R^3}$ such that
\begin{equation}  \label{eq:twotwo}
a+ib=\psi^t\sigma\psi,\; c=-\hat{\psi}^t\sigma\psi.
\end{equation}
Here, the superscript $t$ denotes transposition and $\sigma=(\sigma_1,\sigma_2,\sigma_3)$ is a vector whose cartesian components are the complex symmetric $2\times2$ matrices
\begin{equation}  \label{eq:twothree}
\sigma_1=\left(
           \begin{array}{cc}
             1 & 0 \\
             0 & -1 \\
           \end{array}
         \right),
\sigma_2=\left(
           \begin{array}{cc}
             i & 0 \\
             0 & i \\
           \end{array}
         \right),
\sigma_3=\left(
           \begin{array}{cc}
             0 & -1 \\
             -1 & 0 \\
           \end{array}
         \right).
\end{equation}
In addition to this, $\hat{\psi}$ is the mate (or conjugate) of $\psi$ and $\bar{\psi}$ is complex conjugation of $\psi$ \cite{p}. Therefore,
\begin{equation}  \label{eq:twofour}
\hat{\psi}\equiv-\left(
                  \begin{array}{cc}
                    0 & 1 \\
                    -1 & 0 \\
                  \end{array}
                \right)\cdot \bar{\psi}=\left(
                                        \begin{array}{c}
                                          -\bar{\psi}_2 \\
                                          \bar{\psi}_1 \\
                                        \end{array}
                                      \right).
\end{equation}
In that case, the vectors $a,b$ and $c$ are explicitly given by

\begin{align} \label{eq:twofive}
&a+ib=\psi^t\sigma\psi=(\psi^2_1-\psi^2_2,i(\psi^2_1+\psi^2_2),-2\psi_1\psi_2), \nonumber \\
&c=-\hat{\psi}^t\sigma\psi=(\psi_1\bar{\psi}_2+\bar{\psi}_1\psi_2,i(\psi_1\bar{\psi}_2-\bar{\psi}_1\psi_2),|\psi_1|^2-|\psi_2|^2).
\end{align}
Since the vector $a+ib\in{C}^3$ is an isotropic vector, by means of an easy computation one find that $a,b$ and $c$ are mutually orthogonal \cite{c,cb,bl}. Also, $|a|=|b|=|c|=\bar{\psi}^t\psi$ and $<a\wedge{b},c>$=det$(a,b,c)>0$. On the contrary, for the vectors $a,b$ and $c$ mutually orthogonal vectors of same magnitudes (det$(a,b,c)>0$) there is a spinor defined up to sign such that the equation (\ref{eq:twotwo}) holds.
Under the conditions state above, for two arbitrary spinors $\phi$ and $\psi$, there exist following equalities

\begin{align} \label{eq:twosix}
\overline{\phi^t\sigma\psi}&=-\hat{\phi}^t\sigma\hat{\psi}, \nonumber \\
\widehat{a\phi+b\psi}&=\bar{a}\hat{\phi}+\bar{b}\hat{\psi},
\end{align}
and
\begin{equation*}
\hat{\hat{\psi}}=-\psi.
\end{equation*}
where a and b are complex numbers. The correspondence between spinors and orthogonal bases given by \eqref{eq:twotwo} is two to one; the spinors $\psi$ and $-\psi$ correspond to the same ordered orthonormal bases $\{a,b,c\}$, with $|a|=|b|=|c|$ and $<a\wedge{b},c>$. In addition to that, the ordered triads \emph{\{a,b,c\}}, \emph{\{b,a,c\}} and \emph{\{c,a,b\}} correspond to different spinors.
 Since the matrices $\sigma$ (given by \eqref{eq:twothree}) are symmetric, any pair of spinors $\phi$ and $\psi$ satisfying $\phi^t\sigma\psi=\psi^t\sigma\phi$. The set ${\psi,\hat{\psi}}$  is linearly independent for the spinor $\psi\neq0$ \cite{cb}.

\section{Spinor Darboux Equations}
In this section we investigate the spinor equation of the Darboux equations for a curve on the oriented surface in Euclidean three 3-space $E^3$.

Let \emph{M} be an oriented surface in Euclidean 3-space and let consider a curve $\alpha(s)$ on the surface M. Since the curve $\alpha(s)$ is also in space, there exists Frenet frame \emph{{T,N,B}} at each points of the curve where \emph{T} is unit tangent vector, \emph{N} is principal normal vector and \emph{B} is binormal vector, respectively. The Frenet equations of curve $\alpha(s)$ is given by
\begin{align}\label{eq:threeone}
\frac{dT}{ds}&=\kappa{N} \nonumber\\
\frac{dN}{ds}&=-\kappa{T}+\tau{B} \\
\frac{dB}{ds}&=-\tau{N}\nonumber
\end{align}
where $\kappa$ and $\tau$  are curvature and torsion of the curve $\alpha(s)$, respectively \cite{c}. Unless otherwise stated we assume that the curve  $\alpha(s)$ is a regular curve with arc length parameter. According to the result concerned with the spinor (given by section 2) there exists a spinor $\psi$ such that
\begin{equation}\label{eq:threetwo}
N+iB=\psi^t\sigma\psi,\;T=-\hat{\psi}^t\sigma\psi
\end{equation}
with $\bar{\psi}^t\psi=1$. Thus, the spinor  $\psi$ represent the triad \emph{\{N,B,T\}} and the variations of this triad along the curve $\alpha(s)$  must correspond to some expression for  $\frac{d\psi}{ds}$. That is, the Frenet equations are equivalent to the single spinor equation
\begin{equation}\label{eq:threethree}
\frac{d\psi}{ds}=\frac{1}{2}(-i\tau\psi+\kappa\hat{\psi})
\end{equation}
where $\kappa$ and $\tau$  denote the torsion and curvature of the curve, respectively. The equation \eqref{eq:threethree} is called spinor Frenet equation, \cite{cb}.

Since the curve $\alpha(s)$  lies on the surface \emph{M}, there exists another frame of the curve  $\alpha(s)$ which is called Darboux frame and denoted by \{\emph{T,g,n}\}.
In this frame, \emph{T} is the unit tangent of the curve, \emph{n} is the unit normal of the surface \emph{M} and \emph{g} is a unit vector given by $g=n\wedge{T}$. Since the unit tangent \emph{T} is common in both Frenet frame and Darboux frame, the vectors \emph{N, B, g}, and \emph{n} lie on the same plane. So that, the relations between these frames can be given as follows
\begin{equation}\label{eq:threefour}
\left[
  \begin{array}{c}
    T \\
    g \\
    n \\
  \end{array}
\right]=\left[
          \begin{array}{ccc}
            1 & 0 & 0 \\
            0 & \cos\theta & \sin\theta \\
            0 & -\sin\theta & \cos\theta \\
          \end{array}
        \right]\left[
                 \begin{array}{c}
                   T \\
                   N \\
                   B \\
                 \end{array}
               \right]
\end{equation}
where $\theta$  is the angle between the vectors \emph{g} and \emph{N}. The derivative formulae of the Darboux frame is
\begin{equation}\label{eq:threefive}
\left[
  \begin{array}{c}
    \frac{dT}{ds} \\
    \frac{dg}{ds} \\
    \frac{dn}{ds} \\
  \end{array}
\right]=\left[
          \begin{array}{ccc}
            0 & \kappa_g & \kappa_n \\
            -\kappa_g & 0 & \tau_g \\
            -\kappa_n & -\tau_g & 0 \\
          \end{array}
        \right]\left[
                 \begin{array}{c}
                   T \\
                   g \\
                   n \\
                 \end{array}
               \right]
\end{equation}
where $\kappa_g$ is the geodesic curvature, $k_n$ is the normal curvature and  $\tau_g$ is the geodesic torsion of  $\alpha(s)$ \cite{m}.

Considering the equations \eqref{eq:twoone}, \eqref{eq:twotwo} and \eqref{eq:threetwo} there exists a spinor  $\phi$, defined up to sign, such that
\begin{equation}\label{eq:threesix}
g+in=\phi^t\sigma\phi,\;\;T=-\hat{\phi}^t\sigma\phi
\end{equation}
and
\begin{equation*}
\bar{\phi}^t\phi=1.
\end{equation*}
That is, the spinor $\phi$  represent the triad \{\emph{g,n,T}\} of the curve $\alpha(s)$  on the surface \emph{M}. The variations of the triad \{\emph{g,n,T}\} along curve must correspond to some expression for $\frac{d\phi}{ds}$. Since $\{\phi,\hat{\phi}\}$ is a basis for the two component spinors  $(\phi\neq0)$, there are two functions \emph{f} and \emph{g}, such that
\begin{equation}\label{eq:threeseven}
\frac{d\phi}{ds}=f\phi+g\hat{\phi}
\end{equation}
where the functions f and g are possibly complex-valued functions.
Differentiating the first equation in \eqref{eq:threesix} and using the equation \eqref{eq:threeseven} we have
\begin{align*}
\frac{dg}{ds}+i\frac{dn}{ds}=&\left(\frac{d\phi}{ds}\right)^t\sigma\phi+\phi\sigma\left(\frac{d\phi}{ds}\right)\\
                            =&2f(\phi^t\sigma\phi)+2g(\hat{\phi}^t\sigma\phi).
\end{align*}
Substituting the equations \eqref{eq:threefive} and \eqref{eq:threesix} into the last equation and after simplifying, one finds
\begin{equation}\label{eq:threeeight}
-i\tau_g(g+in)+(-\kappa_g-i\kappa_n)=2f(g+in)-2g(t).
\end{equation}
From the last equation we get
\begin{equation}\label{eq:threenine}
f=-i\frac{\tau_g}{2},\;g=\frac{i\kappa_n+\kappa_g}{2}
\end{equation}
Thus, we have proved the following theorem.
\begin{theorem}
Let the two-component spinor $\phi$  represents the triad
\{g,n,T\} of a curve parametrized by arc length on the
surface M. Then, the Darboux equations are equivalent to the
single spinor equation
\begin{equation}\label{eq:threeten}
\frac{d\phi}{ds}=\left(-i\frac{\tau_g}{2}\right)\phi+\left(\frac{i\kappa_n+\kappa_g}{2}\right)\hat{\phi}
\end{equation}
where $\kappa_n$ and $\kappa_g$  are normal and geodesic
curvature, respectively and $\tau_g$ is the geodesic torsion of
$\alpha(s)$.
\end{theorem}
The equation \eqref{eq:threeten} is called spinor
Darboux equation.
Now, we investigate the relation between the spinors $\psi$ and $\phi$
represent the triad \{\emph{N,B,T}\} (Frenet frame) and
\{\emph{g,n,T}\} (Darboux Frame), respectively.

Considering the equation \eqref{eq:threefour} we get
\begin{equation}\label{eq:threeeleven}
N+iB=(g+in)(\cos\theta+i\sin\theta)
\end{equation}
Thus, from the equations \eqref{eq:threetwo},\eqref{eq:threesix}
and \eqref{eq:threeeleven} we can give following theorem.
\begin{theorem}
Let  the spinors  $\psi$ and $\phi$  represent the triad \{N,B,T\}
and \{g,n,T\} of a curve parametrized with arc length, respectively.
Then, the relation between the spinor formulation of Frenet frame
and Darboux frame is as follows
\begin{align*}
\psi^t\sigma\psi=&e^{i\theta}(\phi^t\sigma\phi) \\
T=&T
\end{align*}
\end{theorem}

\end{document}